\documentclass[runningheads]{llncs}
\usepackage{amsmath}
\usepackage{amsfonts}
\usepackage{mathtools}
\usepackage[T1]{fontenc}
\usepackage{graphicx}
\usepackage[
	pdfstartview=FitH,
	pdffitwindow=true,
	colorlinks,
	linkcolor=black,
	anchorcolor=black,
	citecolor=black,
	urlcolor=blue,
]{hyperref}
\hypersetup
{
	pdftitle     = {Bayesian Line Search},
	pdfsubject   = {University of Hamburg / Bayesian Line Search},
	pdfauthor    = {Robin Labryga},
	pdfkeywords  = {University of Hamburg / Bayesian Line Search},
	plainpages   = false
}
\usepackage{color}

\usepackage{xcolor}
\definecolor{cbBlue}{RGB}{55,  126, 184}
\definecolor{cbOrange}{RGB}{255, 127, 0}
\definecolor{cbGreen}{RGB}{77,  175, 74}
\definecolor{cbPink}{RGB}{247, 129, 191}
\definecolor{cbBrown}{RGB}{166, 86,  40}
\definecolor{cbPurple}{RGB}{152, 78,  163}
\definecolor{cbGray}{RGB}{153, 153, 153}
\definecolor{cbRed}{RGB}{228, 26,  28}
\definecolor{cbYellow}{RGB}{222, 222, 0}
\usepackage[
]{todonotes}
\todostyle{red}{color=red}
\usepackage{tikz}
\usetikzlibrary{positioning}
\usetikzlibrary{calc}
\usepackage{pgfplots}
\usepackage{pgfplotstable}
\usepgfplotslibrary{statistics}
\usepgfplotslibrary{units}
\usepgfplotslibrary{groupplots}
\pgfplotsset{compat=1.18}
\usepackage[capitalise]{cleveref}
\usepackage{algorithm}
\usepackage{algpseudocode}
\usepackage{hhline}
\usepackage{booktabs}

\DeclareMathOperator{\diag}{diag}
\DeclareMathOperator{\EX}{\mathbb{E}}
\DeclareMathOperator{\sign}{sign}
   
\DeclareMathOperator*{\argmax}{argmax}
\newcommand{\ignore}[1]{}
\newcommand{\setR}{\mathbb{R}}

\begin{document}
\title{Information Preserving Line Search via Bayesian Optimization}
\titlerunning{Bayesian Line Search}
%
\author{
	Robin Labryga\inst{1} \and
	Tomislav Prusina\inst{1} \and
	Sören Laue\inst{1}
}
\authorrunning{R. Labryga et al.}
%
\institute{
	University of Hamburg, Germany\\
	\email{\{robin.labryga,tomislav.prusina,soeren.laue\}@uni-hamburg.de}
}
\maketitle              
\begin{abstract}
	Line search is a fundamental part of iterative optimization methods for unconstrained and bound-constrained optimization problems to determine suitable step lengths that provide sufficient improvement in each iteration. Traditional line search methods are based on iterative interval refinement, where valuable information about function value and gradient is discarded in each iteration. We propose a line search method via Bayesian optimization, preserving and utilizing otherwise discarded information to improve step-length choices. Our approach is guaranteed to converge and shows superior performance compared to state-of-the-art methods based on empirical tests on the challenging unconstrained and bound-constrained optimization problems from the CUTEst test set.

	\keywords{Nonlinear Optimization \and Line Search \and Regression \and Bayesian Optimization \and Gaussian Process.}
\end{abstract}
\section{Introduction}
Line search is an essential part of iterative optimization algorithms~\cite{GiesenL19,NW06,sra2012}, determining step lengths that lead to sufficient improvement in each iteration. Traditional methods rely on iterative interval refinement, where an interval containing a potential solution is progressively narrowed by evaluating function values and gradients. While this process systematically reduces the search space, it discards previously obtained information after each iteration. This can result in missing better steps, especially when the refined interval no longer includes a step that offers greater improvement.

This paper introduces a novel approach to inexact line search that leverages Bayesian optimization. Unlike traditional methods, which discard information after each refinement, our approach preserves and utilizes all observed data to make more informed decisions about step lengths. By approximating the objective function with a Gaussian process surrogate model, our method effectively identifies steps that can provide better improvements. It is well-suited for unconstrained or bound-constrained, continuously differentiable black-box objective functions that are bounded from below. Furthermore, the approach is especially advantageous when function and gradient evaluations are expensive.

\subsection{Related Work}

Numerous line search methods have been developed to ensure convergence and improve efficiency. Some widely used methods include backtracking line search~\cite{NW06}, which iteratively reduces the step size until the sufficient decrease condition is met. Hager’s line search~\cite{H89} combines bracketing, bisection, and Newton iterations based on cubic interpolation. The Moré-Thuente line search~\cite{MT94} employs interpolation and an auxiliary function to satisfy the strong Wolfe conditions. The filter line search method by Wächter and Biegler~\cite{WB06} is used in interior-point methods to guide the search by considering both objective function reduction and constraint violations. These methods are widely adopted due to their simplicity and ability to satisfy Wolfe or strong Wolfe conditions~\cite{NW06,W69,W71}. However, they rely on interval refinement, which inherently discards valuable information from intermediate evaluations.

The Barzilai-Borwein line search~\cite{BB88} is a non-monotone step size selection method that estimates the step length using differences in gradients and iterates, inspired by quasi-Newton updates. Unlike traditional line search methods that enforce Wolfe conditions, it adaptively selects step sizes without backtracking, making it particularly effective for quadratic and large-scale unconstrained problems. However, it often requires stabilization to maintain robustness for highly nonlinear objectives.

More recently, Prusina and Laue~\cite{PL24} introduced a surrogate-based line search using cubic interpolation, which retains more information than classical methods. However, this approach remains limited to predefined interpolation models and lacks theoretical convergence guarantees. Similarly, Papageorgiou et al.~\cite{PKK23} proposed a surrogate-based line search method tailored for smooth and derivative-free optimization problems, further expanding the scope of surrogate-based approaches.

Bayesian optimization has emerged as a powerful tool in global optimization, particularly for expensive black-box objective functions~\cite{F18,F18-alt}. It has been successfully applied in various domains, including hyperparameter tuning for machine learning models~\cite{S12}, engineering system design~\cite{F08}, and experimental design selection~\cite{BCF10}. While its use in global optimization is well-documented (e.g., \cite{M94,MM91}), its application to line search processes has received limited attention. Mahsereci and Henning~\cite{MH17} proposed a Bayesian optimization-based line search for stochastic optimization, highlighting the potential of Bayesian methods in this context.

Building on this foundation, our approach leverages Bayesian optimization in the deterministic optimization setting, to systematically model function values and gradients. Unlike classical line search methods that discard intermediate evaluations, our method preserves and utilizes all available information, addressing a critical gap in the literature.

\subsection{Contributions}

This paper makes the following contributions:
\begin{itemize}
    \item \textbf{Novel Line Search Framework}: We propose a Bayesian optimization-based line search method that preserves and utilizes all available data to make more informed step-length decisions.
    \item \textbf{Convergence Proof}: We provide a theoretical proof of convergence for our method, ensuring it can reliably identify sufficient improvements in the optimization process.
    \item \textbf{Extensive Empirical Evaluation}: Our method is integrated into the GENO solver~\cite{LMG19,LaueBG22} and benchmarked against state-of-the-art methods, demonstrating superior performance on challenging unconstrained and bound-con\-strained problems from the CUTEst test set~\cite{GOT15}.
\end{itemize}
The code for our Bayesian optimization-based line search method is available on GitHub\footnote{\url{https://github.com/RobinLabryga/bayesian-geno/tree/bayesian-ls-paper}}.

This paper is organized as follows: \cref{sec:background} provides an overview of line search, Bayesian optimization, and Gaussian processes. \cref{sec:algorithm} details our proposed approach, including its theoretical foundation, convergence proof, and implementation. \cref{sec:experiments} presents empirical evaluations using the GENO solver~\cite{LMG19} on benchmark problems from the CUTEst test set~\cite{GOT15}. Our results demonstrate superior performance compared to state-of-the-art methods.
\section{Background}\label{sec:background}

This section provides an overview of the foundational concepts underlying our approach, including line search, Bayesian optimization, and Gaussian processes. These concepts form the basis for the proposed approach of utilizing Bayesian optimization for line search.

\subsection{Line Search}

Consider a continuously differentiable objective function $f\colon\setR^n \to \setR$ to be minimized.
In iterative optimization, line search starts from a point $x \in \setR^n$ and moves in a descent direction $p$ with $p^\top \nabla f(x) < 0$.
The aim is to find a step length $\alpha$ along $p$ that sufficiently improves upon $f$.
The univariate line search objective $\phi$ and its derivative $\phi'$ along $p$ are defined as
\begin{equation}\label{eq:line_search_objective}
    \phi(\alpha) \coloneq f(x+\alpha\cdot p), \quad \phi'(\alpha) = p^\top \nabla f(x+\alpha\cdot p), \quad \alpha > 0
\end{equation}
Bounds on the step $\alpha$, such as a maximal step $\alpha_{max}$, may also be imposed.

Different line search methods have varying requirements for the chosen step.
The sufficient decrease condition \eqref{eq:sufficient_decrease_condition}, also known as the Armijo rule~\cite{A66}, requires that the function value at $\alpha$ is below a flattened tangent from the starting point, ensuring sufficient decrease:
\begin{equation}\label{eq:sufficient_decrease_condition}
    \phi(\alpha) \leq \phi(0) + \mu \cdot \phi'(0) \cdot \alpha, \quad \mu \in (0,1)
\end{equation}
The curvature condition \eqref{eq:curvature_condition} ensures that the negative gradient at $\alpha$ is not too large compared to the starting point, preventing the step from being too short:
\begin{equation}\label{eq:curvature_condition}
    -\phi'(\alpha) \leq \eta \cdot -\phi'(0), \quad \eta \in (0,1)
\end{equation}
However, the curvature condition does not guarantee that the gradient at $\alpha$ is closer to zero than at the starting point.
The modified curvature condition \eqref{eq:modified_curvature_condition} addresses this by requiring the absolute value of the gradient at $\alpha$ to be sufficiently smaller than at the starting point:
\begin{equation}\label{eq:modified_curvature_condition}
    |\phi'(\alpha)| \leq \eta \cdot |\phi'(0)|, \quad \eta \in (0,1)
\end{equation}
The combination of the sufficient decrease condition \eqref{eq:sufficient_decrease_condition} and the curvature condition \eqref{eq:curvature_condition} is known as the Wolfe conditions, while the combination of the sufficient decrease condition \eqref{eq:sufficient_decrease_condition} with the modified curvature condition \eqref{eq:modified_curvature_condition} is known as the strong Wolfe conditions~\cite{NW06,W69,W71}.
Typically, the parameters $\eta$ and $\mu$ for the (strong) Wolfe conditions are chosen such that $\mu < \eta$ (e.g., $\mu=10^{-4}, \eta=0.9$)~\cite{NW06}.

Optimizers using a line search method that finds a step satisfying the strong Wolfe conditions are convergent if the search direction $p$ is a descent direction.
Gradient methods and quasi-Newton methods satisfy this criterion~\cite{DS96,NW06}.
Additionally, the strong Wolfe conditions ensure a positive definite quasi-Newton update to the Hessian approximation is possible~\cite{NW06}.

To guarantee that a step satisfying the strong Wolfe conditions exists within an interval, the auxiliary function $\psi$~\cite{MT94} is often used:
\begin{equation}\label{eq:auxiliary}
    \psi(\alpha) \coloneq \phi(\alpha) - \left(\phi\left(0\right) + \mu \cdot \phi'\left(0\right) \cdot \alpha\right),
    \qquad \psi'(\alpha) = \phi'(\alpha) - \mu \cdot \phi'(0),
\end{equation}
where $\mu$ is the parameter from the sufficient decrease condition \eqref{eq:sufficient_decrease_condition}.
\ignore{

}
We have the following theorem:
\begin{theorem}[Mor\'{e}-Thuente 2.1~\cite{MT94}]\label{theorem:MT2.1}
    Let $I$ be a closed interval with endpoints $\alpha_l$ and $\alpha_u$.
    If the endpoints satisfy
    \begin{equation*}
        \psi(\alpha_l) \leq \psi(\alpha_u),
        \quad \psi(\alpha_l) \leq 0,
        \quad \psi'(\alpha_l)(\alpha_u - \alpha_l) < 0,
    \end{equation*}
    then there is an $\alpha^\star$ in $I$ with $\psi(\alpha^\star) \leq \psi(\alpha_l)$ and $\psi'(\alpha^\star) = 0$.
    $\alpha^\star$ satisfies the strong Wolfe conditions.
\end{theorem}
From this, we obtain the following corollary:
\begin{corollary}\label{corollary:MT2.1:leftassumptions}
    Let $I=[\alpha_l, \alpha_u]$ be a closed interval with $\alpha_l < \alpha_u$.
    If $\alpha_l$ satisfies
    \begin{equation*}
        \quad \psi(\alpha_l) \leq 0,
        \quad \psi'(\alpha_l) < 0,
    \end{equation*}
    and for $\alpha_u$ we have $\psi(\alpha_u) \geq \psi(\alpha_l)$ or $\psi'(\alpha_u) \geq 0$
    then there is an $\alpha^\star$ in $I$ with $\psi(\alpha^\star) \leq \psi(\alpha_l)$ and $\psi'(\alpha^\star) = 0$.
    $\alpha^\star$ satisfies the strong Wolfe conditions.
\end{corollary}
\begin{proof}
    By \cref{theorem:MT2.1}, if $\psi(\alpha_u)\geq\psi(\alpha_l)$, an $\alpha^\star$ exists in $I$ with $\psi(\alpha^\star) \leq \psi(\alpha_l)$ and $\psi'(\alpha^\star)=0$.
    If $\psi(\alpha_u) < \psi(\alpha_l)$ and $\psi'(\alpha_u)\ge0$, we rename $\alpha_l$ and $\alpha_u$ to complete the argument similarly.
    Thus, $\alpha^\star$ satisfies the strong Wolfe conditions.
    \qed
\end{proof}

\subsection{Bayesian Optimization}\label{sec:bayesian_optimization}
We utilize Bayesian optimization to minimize the line search objective and determine a suitable step length. By approximating the objective function with a Bayesian surrogate, we iteratively refine our model using all information from previously queried points.

Surrogate methods approximate an expensive-to-evaluate objective function by conditioning a surrogate model on known function evaluations. The surrogate model is then used to define an acquisition function that selects the next evaluation point. Since the surrogate is cheaper to evaluate than the objective function, this approach allows for more evaluations when optimizing the acquisition function, ultimately reducing the number of expensive objective function evaluations~\cite{BCF10,F18,F18-alt}. An algorithmic overview of Bayesian optimization is presented in \cref{alg:bayesian-optimization}.

\begin{algorithm}
   \caption{Bayesian Optimization}\label{alg:bayesian-optimization}
   \begin{algorithmic}[1]
      \For{$t\in[1..T]$}
      \State Condition surrogate model $\mathcal{M}$ on $\mathcal{D}_{1:t-1}$
      \State $x_t \gets \argmax_{x \in A} u(x|\mathcal{M})$
      \State $\mathcal{D}_{1:t} \gets \mathcal{D}_{1:t-1} \cup \{x_t\}$
      \EndFor
      \State \Return $\argmax_{x\in\mathcal{D}} \phi(x)$ or $\argmax_{x \in A} \mathcal{M}(x)$
   \end{algorithmic}
\end{algorithm}

Bayesian optimization is typically used for multivariate optimization, but in the context of line search, we focus on the univariate case. Given an objective function $\phi\colon\mathbb{R} \to \mathbb{R}$, our goal is to minimize $\phi(\alpha)$ over a feasible set $A \subset \mathbb{R}$. The Bayesian optimization framework iteratively selects query points to evaluate, refining the surrogate model until a satisfactory solution is found.

Each iteration of Bayesian optimization incurs computational overhead due to updating the surrogate model and optimizing the acquisition function. Consequently, Bayesian optimization is most effective when function evaluations are expensive, and only a limited number of evaluations are feasible~\cite{F18,F18-alt}.

Bayesian optimization is known to converge under mild assumptions on the surrogate model and the acquisition function. Specifically, as the algorithm produces a dense sequence of observations, the surrogate model increasingly refines its approximation of the objective function, leading to improved estimates of the global minimum~\cite{M94,MM91,VB10}.

Common surrogate models in Bayesian optimization include Gaussian processes~\cite{K51,VB10} and Wiener processes~\cite{K64,S88}. In our approach, we use Gaussian processes as the surrogate model, which we discuss in the next section.

\subsection{Gaussian Process}\label{sec:gaussian_process}
A Gaussian process is a collection of random variables where any finite subset follows a joint Gaussian distribution. It is fully specified by a mean function $\mu$  and a covariance function $k$, which define its properties and structure. Formally, a Gaussian process can be written as:
\begin{align*}
   f(x)   & \sim GP(\mu(x), k(p,q))           \\
   \mu(x) & = \EX [f(x)]                       \\
   k(p,q) & = \EX [(f(p)-\mu(p))(f(q)-\mu(q))]
\end{align*}

Before any observations are made, a Gaussian process prior is chosen to reflect prior beliefs about the function’s behavior. When conditioned on observed function values and gradients, the prior is updated to yield a Gaussian process posterior, that incorporates all available data to refine the function approximation.

Given \( n \) observations at points \( X = \{x_1, \dots, x_n\} \), with corresponding function values \( y = \{y_1, \dots, y_n\} \) and (optionally) gradients \( g = \{\nabla f(x_1), \dots, \nabla f(x_n)\} \), the predictive mean $\overline{\mu}$ and covariance $\overline{\text{cov}}$ of the Gaussian process posterior at a new query point \( x \) are given by:
\begin{align*}
   \overline{\mu}(x) & =
   \mu(x) + K(x,X)\left(K(X,X) + \sigma^2\mathbb{I}\right)^{-1}Y
   \\
   \overline{cov}(x) & = K(x,x)-K(x,X)\left(K(X,X) + \sigma^2\mathbb{I}\right)^{-1}K(X,x)
\end{align*}

where \( K(X,X) \) is the Gram matrix of the covariance function, and \( Y \) represents the observed function values. If gradient information is available, the Gaussian process can be conditioned on both function values and gradients to enhance predictive accuracy.

The covariance function, also known as the kernel, plays a crucial role in determining the behavior of the Gaussian process. Different kernels model different assumptions about the smoothness and variability of the approximated function. In our approach, we use the Matérn kernel~\cite{RW06} with \( \nu = \frac{5}{2} \), which provides a balance between smoothness and adaptability:

\[
k_{\nu=5/2}(p,q) = \left(1 + \frac{\sqrt{5}|p-q|}{l} + \frac{5|p-q|^2}{3l^2} \right) e^{-\frac{\sqrt{5}|p-q|}{l}}.
\]

This kernel is particularly effective for modeling functions with finite differentiability, making it well-suited for line search applications. Gaussian processes with the Matérn kernel provide accurate uncertainty quantification~\cite{fiedler2021practical}, which is essential for making informed decisions in Bayesian optimization.

\subsubsection{Differentiability of Gaussian Processes}  
Since we want to optimize acquisition functions defined based on Gaussian processes, having the derivative of Gaussian processes is invaluable.
As differentiation is a linear operator, the derivative of a Gaussian process is another Gaussian process that we can compute~\cite{laue22,RW06}.

\ignore{
The differentiability of a Gaussian process is determined by the smoothness properties of its covariance function. If the kernel function is infinitely differentiable (such as the squared exponential kernel), the resulting Gaussian process is also infinitely differentiable. However, kernels like the Matérn family introduce flexibility in smoothness. The Matérn kernel with parameter \( \nu \) ensures that sample paths of the Gaussian process are \( \lfloor \nu - 1/2 \rfloor \) times differentiable in the mean-square sense. For the Matérn \( \nu = 5/2 \) kernel used in our approach, this guarantees that the function is at least twice differentiable, making it suitable for optimization methods that require gradient information. The differentiability property of the Gaussian process plays a crucial role in conditioning on gradient observations, enabling more accurate modeling of function behavior while preserving computational efficiency.
}

The derivative of the predictive mean and the derivative of the predictive covariance can be calculated via
\begin{align*}
   \frac{\partial}{\partial x}\overline{\mu}(x)
                                                & = \frac{\partial}{\partial x}\mu(x) + \frac{\partial}{\partial x}K(x,X)\left(K(X,X) + \sigma^2\mathbb{I}\right)^{-1}Y
   \\
   \frac{\partial}{\partial x}\overline{cov}(x) & = \frac{\partial}{\partial x}K(x,x) - 2 \cdot \frac{\partial}{\partial x}K(x,X) \left(K(X,X) + \sigma^2\mathbb{I}\right)^{-1} K(x,X)^\top \\
\end{align*}

\subsection{Acquisition Functions}
In Bayesian optimization, an acquisition function determines the next point to evaluate by balancing exploration (searching in uncertain regions) and exploitation (focusing on regions likely to yield optimal values). This process can also be viewed as minimizing a risk function~\cite{M94}.
The choice of acquisition function significantly affects the efficiency of the optimization process. 

Examples of acquisition functions include Probability of Improvement ($PI$)~\cite{K64}, Expected Improvement ($EI$)~\cite{JSW98}, Lower Confidence Bound ($LCB$)~\cite{CJ92}, knowledge gradient ($KG$)~\cite{FPD09}, derivative-enabled Expected Improvement ($d-EI$)~\cite{WPWF17}, and derivative-enabled knowledge gradient ($d-KG$)~\cite{WPWF17}.
Variants of these acquisition functions, such as $GP-LCB$~\cite{SKKS09}, also exist.

Most of these acquisition functions are designed for derivative-free scenarios, with $d-EI$ and $d-KG$ being the exceptions that explicitly incorporate gradient information.
Since acquisition functions depend on the mean and variance of the surrogate Gaussian process, gradients implicitly affect all acquisition functions if gradient information is used to condition the Gaussian process posterior.

The variance is the diagonal of the covariance, i.e., $\mathbb{V}(x)=\diag(\overline{cov}(x))$, and the standard deviation is the square root of the variance, i.e., $\mathbb{S}(x)=\sqrt{\mathbb{V}(x)}$. Then, the Lower Confidence Bound ($LCB$) acquisition function~\cite{CJ92} and its derivative are defined as follows:
\begin{align*}
    LCB(x)                             & = -\mu(x) + \kappa \mathbb{S}(x)                                                        \\
    \frac{\partial}{\partial x} LCB(x) & = -\frac{\partial}{\partial x}\mu(x) + \kappa \frac{\partial}{\partial x} \mathbb{S}(x)
\end{align*}
where $\kappa\geq0$ is a hyperparameter. A larger value of $\kappa$ increases exploration by prioritizing points with high uncertainty, while a smaller $\kappa$ focuses more on exploitation by selecting points with lower predicted function values. The LCB acquisition function is particularly useful when minimizing the function with confidence bounds. In our approach, we utilize the LCB acquisition function.
\section{Algorithm}\label{sec:algorithm}

\begin{figure}[]
   \centering
   \begin{tikzpicture}
      \begin{groupplot}[
            group style={
                  group size=1 by 2,
                  vertical sep=0cm,
                  xlabels at=edge bottom,
                  xticklabels at=edge bottom,
               },
            ytick align=outside,
            ytick pos=left,
            width=\textwidth,
            xmajorgrids=true,
            xmin=0-.1,
            xmax=1+.1,
         ]

         \nextgroupplot[
            ylabel={$\phi$},
            xlabel={},
            xtick pos=top,
            height=0.5\textwidth,
            legend style={
                  legend cell align={left},
                  legend pos=north west,
               },
         ]
         \addplot[cbBlue, thick, dotted] table [col sep=comma, x=x, y=y]{data/compare_MT_BO/data.csv};
         \addlegendentry{$\phi(\alpha)$};
         \addplot[cbRed, thick] table [col sep=comma, x=x, y=mean,]{data/compare_MT_BO/data.csv};
         \addplot[cbRed, dashed, forget plot] table [col sep=comma, x=x, y=std_upper,]{data/compare_MT_BO/data.csv};
         \addplot[cbRed, dashed, forget plot] table [col sep=comma, x=x, y=std_lower,]{data/compare_MT_BO/data.csv};
         \addlegendentry{$\overline{\mu}(\alpha) \pm 2 \cdot \mathbb{S}(\alpha)$};
         \addplot[only marks, mark=*, cbRed] coordinates {(0,1.0) (1,2.0) (0.31179865369365983,0.29394293171002456) (0.22454876024213324,0.28087988873303443) (0.5841504513200902,0.6146751448527172) (0.7277287189554162,0.0695822941436189)};
         \addlegendentry{Bayesian optimization observations};
         \addplot[only marks, mark=*, cbGreen] coordinates {(0.0,1.0) (1.0,2.0) (0.4040903836973557,0.5725692775683063) (0.21937338528328285,0.28735118604218773) (0.26773397163991036,0.2566737683873239)};
         \addlegendentry{Mor\'{e}-Thuente observations};

         \nextgroupplot[
            ylabel={$LCB$},
            xlabel={$\alpha$},
            height=0.25\textwidth,
            ytick align=outside,
            ytick pos=left,
            xtick pos=bottom,
         ]
         \addplot[cbOrange, thick, forget plot] table [col sep=comma, x=x, y=lcb,]{data/compare_MT_BO/data.csv};

      \end{groupplot}
   \end{tikzpicture}
   \caption{
   An exemplary line search comparing the Mor\'{e}-Thuente line search with Bayesian optimization.
   The Mor\'{e}-Thuente line search fails to find the global optimizer, while Bayesian optimization is on track to locate it.
   The orange line represents the Lower Confidence Bound ($LCB$) acquisition function, defined based on the Gaussian process conditioned on all prior observations.
   We observe that the next observation by Bayesian optimization, according to the acquisition function, is very close to the global optimizer.
   }
   \label{fig:MT-vs-BO}
\end{figure}
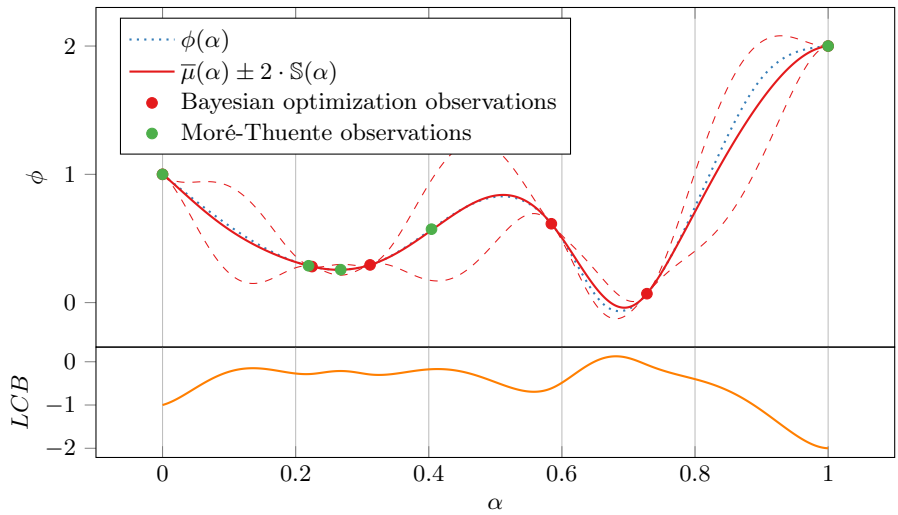

We apply Bayesian optimization to the line search objective \( \phi(\alpha) \) (\cref{eq:line_search_objective}) to find a step size that satisfies the strong Wolfe conditions. Unlike traditional line search methods that discard intermediate evaluations, our approach retains and utilizes all function and gradient information observed during the search. This is achieved by using a Gaussian process model with the Matérn covariance function and selecting step sizes via the Lower Confidence Bound (LCB) acquisition function. An illustrative comparison between the Moré-Thuente line search and Bayesian optimization-based line search is shown in \cref{fig:MT-vs-BO}.

While Bayesian optimization offers a more global approach to finding optimal step sizes, it may suffer from slow convergence if applied naively. Additionally, Bayesian optimization requires an interval with finite endpoints, whereas traditional line searches often operate over unbounded intervals (e.g., \( [0, \infty) \)). We address these challenges by dynamically refining the search interval and integrating Bayesian optimization with an interval update mechanism.

\subsection{Interval Selection}

To apply Bayesian optimization effectively, we first determine an initial bounded interval that contains a step satisfying the strong Wolfe conditions. We initialize the search interval as \( [\alpha_l = 0, \alpha_u = \alpha_0] \), where \( \alpha_0 = \min(1, \alpha_{\max}) \). If this interval does not contain a valid step, we expand it iteratively by updating the interval as follows:
\begin{equation*}
   [\alpha_l, \alpha_u] \leftarrow [\alpha_u, \min(c \cdot \alpha_u, \alpha_{\max})],\quad c \in \mathbb{R}_{>1}
\end{equation*}

If $\psi(\alpha_u) \geq \psi(\alpha_l)$ or $\psi'(\alpha_u) \geq 0$, we can guarantee the presence of a step satisfying the strong Wolfe conditions within the interval according to \cref{corollary:MT2.1:leftassumptions}, since $\alpha_l < \alpha_u$, $\psi(\alpha_l)\leq0$, and $\psi'(\alpha_l) < 0$.
Note that each interval update implies an improvement upon $\phi$, as $\psi(\alpha_u)<\psi(\alpha_l)$ and $\alpha_l < \alpha_u$ implies $\phi(\alpha_u)<\phi(\alpha_l)$.

\subsection{Bayesian Optimization for Step Selection}

Once we have an interval \( [\alpha_l, \alpha_u] \) that either contains a valid step or has reached the maximum step size, i.e., $\alpha_u=\alpha_{\max}$, we apply Bayesian optimization to approximate the minimizer of \( \phi(\alpha) \). This process continues until a step satisfying the strong Wolfe conditions is found or the number of function evaluations exceeds a predefined limit. If Bayesian optimization does not identify a valid step, we refine the interval using a strategy inspired by the Moré-Thuente update rules and retry Bayesian optimization on the new interval.

Unlike traditional methods that consider only two previous evaluations, Baye\-sian optimization tracks all function and gradient values observed during the search. We maintain a dataset \( D \) of all previously evaluated steps. If at any point a step satisfying the strong Wolfe conditions is found and offers better improvement than prior steps, the line search terminates immediately.

During interval refinement, we define an auxiliary function \( \Psi(\alpha) \) to determine the next nested interval. Initially, we set \( \Psi = \psi \), but if we encounter a step \( \alpha_t \) where \( \psi(\alpha_t) \leq 0 \) and \( \phi'(\alpha_t) > 0 \), we switch to \( \Psi = \phi \), following the Moré-Thuente approach.

In each iteration, Bayesian optimization selects the next candidate step \( \alpha_t \) within \( (\alpha_l, \alpha_u) \). If Bayesian optimization fails to improve upon the interval endpoints, we select the step in $\mathcal{D}$ with the highest Gaussian kernel density estimation~\cite{S15} within the current interval, as Bayesian optimization tends to produce denser observations around minimizers asymptotically~\cite{M94,MM91}. Since $\alpha_t$ is now distinct from $\alpha_l$ and $\alpha_u$, we update the interval according to Mor\'{e}-Thuente's updating rules (U1-U3):
\begin{description}
   \item[Case U1:] If $\Psi(\alpha_t) > \Psi(\alpha_l)$, then $\alpha_l^+ \gets \alpha_l$ and $\alpha_u^+ \gets \alpha_t$
   \item[Case U2:] If $\Psi(\alpha_t) \leq \Psi(\alpha_l)$ and $\Psi'(\alpha_t)(\alpha_l - \alpha_t) > 0$, then $\alpha_l^+ \gets \alpha_t$ and $\alpha_u^+ \gets \alpha_u$
   \item[Case U3:] If $\Psi(\alpha_t) \leq \Psi(\alpha_l)$ and $\Psi'(\alpha_t)(\alpha_l - \alpha_t) < 0$, then $\alpha_l^+ \gets \alpha_t$ and $\alpha_u^+ \gets \alpha_l$
\end{description}

This process is repeated until a step satisfying the strong Wolfe conditions is found or the iteration limit is reached. To ensure finite convergence, we impose a condition on the interval size: if the interval has not decreased by a factor of \( \delta < 1 \) (typically \( \delta = 2/3 \)) over two consecutive refinements, we force a bisection step.
Mor\'{e}-Thuente show that this process returns a step satisfying the strong Wolfe conditions or $\alpha_{max}$.

\subsection{Bayesian Optimization Implementation}

The Bayesian optimization approximates the minimizer of $\phi$ on an interval $[\alpha_{min}, \alpha_{max}]$ with $\alpha_{min} < \alpha_{max}$.
Bayesian optimization runs until a set threshold on the number of condition steps is reached or a step satisfying the strong Wolfe conditions is found.
In each iteration, we condition the Gaussian process used as the Bayesian optimization surrogate on all function values and gradients of the steps in $\mathcal{D} \cap [\alpha_{min}, \alpha_{max}]$.
Our Gaussian process prior uses the Mat\'ern kernel with $\nu=\frac{5}{2}$ and length scale hyperparameter $|\alpha_{max}-\alpha_{min}|$.
The prior mean is a constant at the best-known function value ($\mu(x)=\min \{\phi(\alpha)|\alpha\in\mathcal{D} \cap [\alpha_{min}, \alpha_{max}]\}$), as we have no other information about the structure of the black-box objective function.

We maximize the Lower Confidence Bound ($LCB$) acquisition function using a hybrid global-local optimizer that employs the DIRECT~\cite{JM21,JPS93} optimizer to find a step and refines it using L-BFGS-B~\cite{ZBLN97}.
The new step is added to $\mathcal{D}$ and used to condition the Gaussian process in the next iteration.

If the number of steps used to condition the surrogate Gaussian process exceeds a parameter threshold without identifying a step satisfying the strong Wolfe conditions, we return the step in $\mathcal{D} \cap [\alpha_{min}, \alpha_{max}]$ with the best function value.

\subsection{Convergence Guarantee}

To prove that the discovery of the initial interval with finite endpoints terminates in a finite number of iterations, we show that the number of interval updates required until a step satisfying the strong Wolfe conditions can be guaranteed to be within the interval or $\alpha_u=\alpha_{max}$ is bounded from above.
\begin{theorem}
   Let $\phi_{min}$ be a strict lower bound for $\phi$.
   Let $[\alpha_l=0,\alpha_u=\alpha_0=\min(1, \alpha_{max})]$ be the initial interval.
   The number of interval updates of the form $[\alpha_u, \min(c \cdot \alpha_u, \alpha_{max})], c\in\mathbb{R}_{>1}$, until $\alpha_u=\alpha_{max}$ or $\psi(\alpha_u) \geq \psi(\alpha_l)$ or $\psi'(\alpha_u) \geq 0$ is bounded from above by
   \begin{equation*}
      \left\lceil\frac{1}{\alpha_0}\log_c\left(\min\left(\frac{1}{\mu} \cdot \frac{\phi_{min} - \phi(0)}{\phi'(0)}, \alpha_{max}\right)\right)\right\rceil
   \end{equation*}
\end{theorem}
\begin{proof}
   For all $\alpha > \alpha_b$ defined as
   \begin{equation*}
      \alpha_b \coloneq \frac{1}{\mu} \cdot \frac{\phi_{min} - \phi(0)}{\phi'(0)}
   \end{equation*}
   we have
   \begin{align*}
      \psi(\alpha) & = \phi(\alpha) - (\phi(0) + \mu \cdot \phi'(0) \cdot \alpha)                                                             \\
                   & > \phi_{min} - \left(\phi\left(0\right) + \mu \cdot \phi'\left(0\right) \cdot \alpha_b\right)                            \\
                   & = \phi_{min} - \left(\phi(0) + \mu \cdot \phi'(0) \cdot \frac{1}{\mu} \cdot \frac{\phi_{min} - \phi(0)}{\phi'(0)}\right) \\
                   & = \phi_{min} - \phi_{min} = 0
   \end{align*}
   This means that no $\alpha > \alpha_b$ can satisfy the sufficient decrease condition.
   Suppose we update the interval as specified enough, such that $\alpha_u > \alpha_b$, without any prior interval satisfying $\psi(\alpha_u) \geq \psi(\alpha_l)$ or $\psi'(\alpha_u) \geq 0$.
   Then we have $\alpha_u > \alpha_b$, $\psi(\alpha_l) \leq \psi(0) \leq 0$, $\psi'(\alpha_l)<0$ and $\psi(\alpha_u) > 0$.
   Thus, we have found an interval with $\psi(\alpha_u) \geq \psi(\alpha_l)$.
   The number of interval updates required to reach this interval or have $\alpha_u = \alpha_{max}$ is $\left\lceil\frac{1}{\alpha_0}\log_c\left(\min\left(\alpha_b, \alpha_{max}\right)\right)\right\rceil$.
   \qed
\end{proof}

Since the interval found using the above process satisfies the invariant of Mor\'{e}-Thuente up to reordering, the interval refinement process terminates in a finite number of iterations.
Combining the convergence guarantee of the initial interval selection and the convergence guarantee of the interval refinement process, we can ensure that our Bayesian line search terminates after a finite number of iterations, returning either a step satisfying the strong Wolfe conditions or $\alpha_{\max}$.
\section{Experiments}\label{sec:experiments}

To evaluate the performance of our Bayesian optimization-based line search method described in \cref{sec:algorithm}, we integrate it into the GENO solver~\cite{LMG19}, which employs a quasi-Newton optimization approach, and replace its existing line search component.

We conduct experiments similar to those by Prusina and Laue~\cite{PL24}, comparing the performance of our method in terms of convergence, the number of function evaluations, and runtime.
In addition to the GENO solver with Bayesian line search (Bayesian GENO), we evaluate:
\begin{itemize}
    \item The GENO solver with cubic spline interpolation line search from Prusina and Laue~\cite{PL24} (Cubic GENO),
    \item The L-BFGS-B solver~\cite{ZBLN97} interfaced via SciPy~\cite{scipy20}, which uses Moré-Thuente line search~\cite{MT94},
    \item The Ipopt solver, which employs a filter line search~\cite{WB05,WB06}.
\end{itemize}
While Bayesian GENO and Cubic GENO modify only the line search component of the GENO solver, the Ipopt and L-BFGS-B experiments use entirely different solver implementations.
As a result, performance differences may partly reflect solver-specific factors, while still offering useful comparative insights.

The experiments are conducted on 292 unconstrained and 163 bound-constrained problems from the CUTEst test set~\cite{GOT15}, with a time limit of 3000 seconds per problem. Three unconstrained problems with an unbounded optimal solution are excluded.

The experiments are executed on a machine with an Intel Xeon Gold 5315Y CPU (3.20 GHz), 256 GiB of RAM, and Ubuntu 22.04.4 LTS. All solvers access objective function values and gradients, starting from the same initial point \( x_0 \).

The code for our benchmark is available on GitHub\footnote{\url{https://github.com/RobinLabryga/bayesian-geno-benchmark/tree/bayesian-ls-paper}}.

\subsection{Convergence}

A solver is considered converged if it satisfies at least one of the following conditions:

\begin{enumerate}
    \item Function value convergence, defined as:
	\begin{equation}\label{eq:convergence_value}
	    \frac{f_{solver}-f^\star}{1+|f^\star|}<10^{-4}
	\end{equation}
    where \( f_{\text{solver}} \) is the function value obtained by the solver, and \( f^* \) is the best function value found by any solver, serving as an approximation of the global optimum.
    
    \item Gradient convergence, defined as:
	\begin{equation}\label{eq:convergence_gradient}
	    \frac{\|g\|_\infty}{1+|f_{solver}|}<10^{-6}
	\end{equation}
    where \( \|\nabla g\|_{\infty} \) is the infinity norm of the gradient at the final solution point.
\end{enumerate}

To ensure valid comparisons, any solver that violates bound constraints or returns an inconsistent function value is corrected by clipping the solution to the feasible bounds and recomputing the function value and gradient.

\Cref{table:solver-convergence} presents the number of problems successfully solved by each method for both unconstrained and bound-constrained cases. The best values in each column are highlighted.

\begin{table}
    \centering
    \caption{
        The number of unconstrained and bound-constrained problems each method solved.
        Column $f$ conv. shows how many problems each method solved by the function value criterion (\cref{eq:convergence_value}), while column $g$ conv. shows how many problems each solver solved by the gradient criterion (\cref{eq:convergence_gradient}).
        Column conv. shows the number of problems where at least one of the two criteria is satisfied.
        The best value in each column is bold.
    }
    \setlength{\tabcolsep}{3pt}
    \begin{tabular}{l@{\hspace{0.2cm}}rrr@{\hspace{0.5cm}}rrr}
        \toprule
        \multicolumn{1}{c}{} & \multicolumn{3}{c}{unconstrained} & \multicolumn{3}{c}{bound-constrained}                                                             \\
        Solver               & $f$ conv.
                             & $g$ conv.
                             & conv.
                             & $f$ conv.
                             & $g$ conv.
                             & conv                                                                                                                                  \\
        \midrule
        Bayesian GENO        & \textbf{256}                      & \textbf{222}                          & \textbf{269} & \textbf{140} & \textbf{145} & \textbf{156} \\
        Cubic GENO           & 252                               & 213                                   & 263          & 125          & 127          & 142          \\
        Ipopt                & 195                               & 194                                   & 236          & 124          & 93           & 136          \\
        L-BFGS-B             & 232                               & 217                                   & 246          & 135          & 139          & 146          \\
        \bottomrule
    \end{tabular}
    \label{table:solver-convergence}
\end{table}

Our line search method outperforms all other solvers in all three categories.
The Bayesian GENO solver converges on more problems in function value, gradient, or at least one criterion compared to the other solvers for both unconstrained and bound-constrained problems.

\subsection{Function Evaluations}

Since Bayesian optimization promises a lower function evaluation requirement to approximate a minimizer, we compare the number of function evaluations used by each solver.
\Cref{fig:feval-comparison} (a) and (b) compare the solvers regarding the number of function evaluations for unconstrained and bound-constrained problems respectively.

\newcommand{\addfevalplot}[3]{
    \addplot[#3, thick] table [col sep=comma, x=feval, y=mean,]{#1};
    \addlegendentry{#2 $\mu \pm \sigma$}
    \addplot[#3, dashed, forget plot] table [col sep=comma, x=feval, y=std_upper,]{#1};
}

\begin{figure}[]
    \centering
    \begin{minipage}[]{\textwidth}
        \centering
        \begin{tikzpicture}
            \begin{axis}[
                    width=\textwidth,
                    height=.5\textwidth,
                    legend cell align={left},
                    ymode=log,
                    xlabel={Function evaluations},
                    ylabel={Distance to $f^\star$},
                    xmin=0-.1,
                    xmax=4+.1,
                    ymax=1+.025,
                    xtick={0,1,2,3,4},
                ]
                \addfevalplot{data/unconstrained/feval_normalizedBayesianGenoSolver.csv}{Bayesian GENO}{cbRed}
                \addfevalplot{data/unconstrained/feval_normalizedCubicGenoSolver.csv}{Cubic GENO}{cbBlue}
                \addfevalplot{data/unconstrained/feval_normalizedIpoptSolver.csv}{Ipopt}{cbOrange}
                \addfevalplot{data/unconstrained/feval_normalizedScipySolver.csv}{L-BFGS-B}{cbGreen}
            \end{axis}
        \end{tikzpicture}
        \textbf{(a)} Unconstrained problems
    \end{minipage}
    \hfill
    \begin{minipage}[]{\textwidth}
        \centering
        \begin{tikzpicture}
            \begin{axis}[
                    width=\textwidth,
                    height=.5\textwidth,
                    legend cell align={left},
                    ymode=log,
                    xlabel={Function evaluations},
                    ylabel={Distance to $f^\star$},
                    xmin=0-.1,
                    xmax=4+.1,
                    ymax=1+.025,
                    xtick={0,1,2,3,4},
                ]
                \addfevalplot{data/bound/feval_normalizedBayesianGenoSolver.csv}{Bayesian GENO}{cbRed}
                \addfevalplot{data/bound/feval_normalizedCubicGenoSolver.csv}{Cubic GENO}{cbBlue}
                \addfevalplot{data/bound/feval_normalizedIpoptSolver.csv}{Ipopt}{cbOrange}
                \addfevalplot{data/bound/feval_normalizedScipySolver.csv}{L-BFGS-B}{cbGreen}
            \end{axis}
        \end{tikzpicture}
        \textbf{(b)} Bound-constrained problems
    \end{minipage}
    \caption{
        A comparison between the number of function evaluations required by each solver in the experiments.
        For each problem, we calculate the distance of the best-observed function value to $f^\star$ at every function evaluation.
        We then normalize the distance to $f^\star$ such that the distance at the first function evaluation is $1.0$.
        Similarly, we normalize the number of function evaluations such that the number of function evaluations at the discovery of $f^\star$ is $1.0$.
        We then calculate the mean (solid lines) and standard deviation (dashed lines) of the normalized distances across all problems.
        The $x$-axis shows the number of function evaluations, while the $y$-axis shows the distance to $f^\star$.
    }
    \label{fig:feval-comparison}
\end{figure}
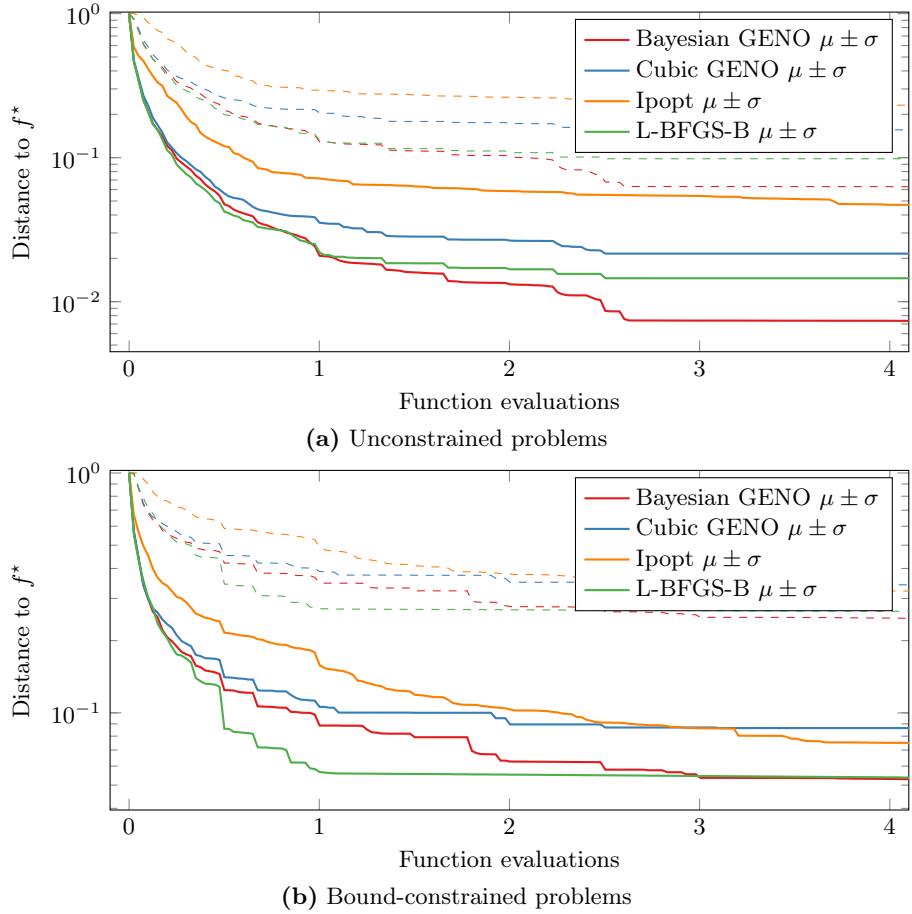

The Bayesian GENO solver requires fewer function evaluations than Ipopt and Cubic GENO to reach the same distance to the optimum $f^\star$.
Compared to L-BFGS-B, the Bayesian GENO solver is initially on par but outperforms L-BFGS-B with increasing function evaluations on unconstrained problems.
On bound-constrained problems, the Bayesian GENO solver initially requires more function evaluations than L-BFGS-B but outperforms it with increasing function evaluations.
This is expected since Bayesian GENO uses a more global exploration approach, requiring more function evaluations per line search.
The improvement over Cubic GENO indicates a more efficient exploration strategy in Bayesian GENO.

\subsection{Overhead}

Each solver incurs different overhead in addition to the time required to evaluate the objective function.
We compare the overhead incurred by each solver.

\begin{figure}[]
    \centering
    \begin{minipage}[]{.45\textwidth}
        \centering
        \begin{tikzpicture}
            \begin{axis}[
                    height=.67\textwidth,
                    width=0.9\textwidth,
                    ytick={1,2,3,4},
                    yticklabels={L-BFGS-B, Ipopt, Cubic GENO, Bayesian GENO},
                    xmin=-1,
                    xmax=16,
                    xtick={0,5,10,15},
                    xlabel={Overhead},
                    x unit=s, x unit prefix=m,
                ]
                \addplot+ [cbGreen, boxplot prepared={
                            lower whisker=0.01340401172637948, lower quartile=0.0594132382448598, median=0.06962907101349, upper quartile=0.48865602611366554, upper whisker=1.1281548775315917
                        }] coordinates {};

                \addplot+ [cbOrange, boxplot prepared={
                            lower whisker=0.0, lower quartile=0.43729851394891733, median=0.7747911815909115, upper quartile=3.9158377749109046, upper whisker=9.019325729459524
                        }] coordinates {};

                \addplot+ [cbBlue, boxplot prepared={
                            lower whisker=0.0, lower quartile=0.1625487849602015, median=0.1897808677916016, upper quartile=0.4323924300382443, upper whisker=0.8245010646594133
                        }] coordinates {};

                \addplot+ [cbRed, boxplot prepared={
                            lower whisker=0.0, lower quartile=0.36362162590887726, median=0.7947969623076406, upper quartile=2.8898598846162113, upper whisker=6.178842591294488
                        }] coordinates {};
            \end{axis}
        \end{tikzpicture}
        \textbf{(a)} Unconstrained problems
    \end{minipage}
    \hfill
    \begin{minipage}[]{.45\textwidth}
        \centering
        \begin{tikzpicture}
            \begin{axis}[
                    height=.67\textwidth,
                    width=0.9\textwidth,
                    ytick={1,2,3,4},
                    yticklabels={L-BFGS-B, Ipopt, Cubic GENO, Bayesian GENO},
                    xmax=16,
                    xmin=-1,
                    xtick={0,5,10,15},
                    yticklabel pos=right,
                    xlabel={Overhead},
                    x unit=s, x unit prefix=m,
                ]
                \addplot+ [cbGreen, boxplot prepared={
                            lower whisker=0.0, lower quartile=0.05905301035153451, median=0.0755374149001878, upper quartile=0.6357999521986211, upper whisker=1.2237113965306916
                        },
                ] coordinates {};

                \addplot+ [cbOrange, boxplot prepared={
                            lower whisker=0.0, lower quartile=0.6214537540763891, median=0.9861927076496861, upper quartile=6.208377536179906, upper whisker=14.245551875589324
                        }] coordinates {};

                \addplot+ [cbBlue, boxplot prepared={
                            lower whisker=0.030454523861408234, lower quartile=0.16428654607046736, median=0.19882636479273136, upper quartile=0.4978434079920989, upper whisker=0.9876312031382068
                        }] coordinates {};

                \addplot+ [cbRed, boxplot prepared={
                            lower whisker=0.0, lower quartile=0.7224261716350415, median=2.0258221895330486, upper quartile=6.5033243989625635, upper whisker=14.953257427831373
                        }] coordinates {};
            \end{axis}
        \end{tikzpicture}
        \textbf{(b)} Bound-constrained problems
    \end{minipage}
    \caption{
        Box plots depicting the overhead per function evaluation in milliseconds for each solver on the unconstrained and bound-constrained problems computed as in \cref{eq:overhead}.
    }
    \label{fig:solver-times}
\end{figure}
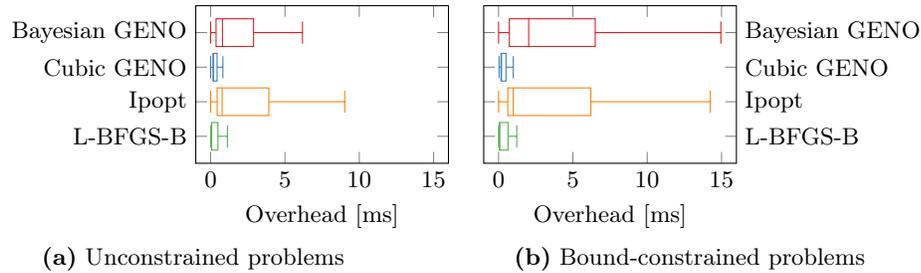

The overhead per function evaluation of a solver for a problem is calculated as
\begin{equation}\label{eq:overhead}
    \frac{t_{\text{solver}}}{n_{\text{solver}}} - \mathbb{E}[t_{\text{problem}}]
\end{equation}
where $t_{\text{solver}}$ is the total time the solver took to solve the problem, $n_{\text{solver}}$ is the number of function evaluations the solver performed, and $\mathbb{E}[t_{\text{problem}}]$ is the expected time to evaluate the objective function once.
The expected time $\mathbb{E}[t_{\text{problem}}]$ is estimated by measuring the runtime of 1000 function evaluations and taking the average.
The box plots in \cref{fig:solver-times} depict the range of overhead estimations across the problems for each solver.

The Bayesian GENO solver has a higher overhead per function evaluation than the Cubic GENO and L-BFGS-B solvers.
The overhead of Bayesian GENO is comparable to that of Ipopt.
The added complexity of Gaussian process regression compared to cubic spline interpolation may explain the performance difference between Bayesian GENO and Cubic GENO.
This is supported by the increased overhead in the bound-constrained setting compared to the unconstrained setting, where Bayesian optimization may search for a non-existent step satisfying the strong Wolfe conditions whenever the search interval includes $\alpha_{max}$.

\section{Conclusion}

We have introduced a novel line search method that uses  Bayesian optimization to better utilize information from function and gradient evaluations. Our approach addresses the limitations of traditional line search techniques, which discard valuable information during iterative refinement. By employing a Gaussian process surrogate model, our method finds more informed step lengths, offering a solution that is both theoretically sound and practically efficient.

Through extensive empirical evaluations on the challenging unconstrained and bound-constrained optimization problems from the CUTEst test set, we demonstrated the superiority of the proposed method when integrated into the GENO solver. The Bayesian line search outperformed state-of-the-art methods in terms of convergence on both function value and gradient criteria, solving a greater number of problems across diverse problem types. While it incurs higher computational overhead, this method is especially valuable in optimization scenarios where function evaluations are expensive.

\begin{credits}

	\subsubsection{\discintname}
	The authors have no competing interests to declare that are relevant to the content of this article.
\end{credits}
%
%

\bibliographystyle{splncs04}
\bibliography{bibliography}

\end{document}